\documentclass[a4paper, oneside, notitlepage, 10pt]{amsart}

\usepackage{amsmath}
\usepackage{amsfonts}
\usepackage{amssymb}
\usepackage{amsthm}
\usepackage{graphicx}
\usepackage{mathrsfs}

\usepackage{verbatim}
\usepackage[T1]{fontenc}        
\usepackage[english]{babel}                              
\usepackage[latin1]{inputenc}  

\newcommand{\sra}{\rightarrow}             

\renewcommand{\l}{\langle}
\renewcommand{\r}{\rangle}

\newcommand{\N}{\mathbb{N}}  
\newcommand{\Z}{\mathbb{Z}}   
\newcommand{\R}{\mathbb{R}}  

\newcommand{\mas}[1]{\left\{ #1 \right\}}
\newcommand{\pan}[1]{\left(  #1 \right)}
\newcommand{\hak}[1]{\left[  #1 \right]}

\newcommand{\phib}{\varphi}

\newcommand{\T}{\mathbb{T}}

\renewcommand{\l}{\langle}
\renewcommand{\r}{\rangle}

\theoremstyle{plain} 
\newtheorem{thm}{Theorem}
\newtheorem{lem}[thm]{Lemma}

\theoremstyle{definition}

\theoremstyle{remark} 
\newtheorem*{rema}{Remark}

\newtheorem{exem}{Example}

\begin{document}
\title{Continuous Measures on Homogenous Spaces}
\author{Michael Bj\"orklund and Alexander Fish}
\date{}

\begin{abstract}
In this paper we generalize Wiener's characterization of continuous measures to compact homogenous manifolds. In particular, we give necessary and sufficient conditions on 
probability measures on compact semisimple Lie groups and nilmanifolds to be continuous. The methods use only simple properties of heat kernels. 
\end{abstract}

\maketitle

\section{Introduction}

N. Wiener gave a very simple, necessary and sufficient, condition for the support of a probability measure $ \mu $ on the unit circle $ \T $ to contain points with positive mass. In fact, he proves the following formula (\cite{GM79})
\[
\lim_{N \sra \infty} \frac{1}{2N+1} \sum_{|k| \leq N} |\hat{\mu}(k)|^2 = \sum_{x \in \T} |\mu(\mas{x})|^2.
\]
This result has found many applications in various areas in mathematics; most notably in ergodic theory and optimal control theory. 

In this note we extend Wiener's criterion to compact homogenous Riemannian manifolds, i.e. compact Riemannian manifolds which admit  a transitive action by  isometries. Our two main
examples are symmetric spaces of compact type ( where other extensions of Wiener's lemma are known, due to Anoussis and \mbox{Bisbas \cite{AnBi00} )} and nilmanifolds. The Fourier coefficient of 
the measure is replaced with the integral of an eigenfunction of the Laplacian against the given measure. We will make use of the short--time behaviour of heat kernels on compact manifolds to deduce 
our results. 

\section{A Generalized Wiener Lemma}

\subsection{Approximations of Diagonal Measures}

In this section we will briefly discuss the simple idea which lies behind this paper. Let $ ( X , \mu ) $ be a probability space, 
and suppose there is a family of complex--valued measurable functions on the direct product $ X \times X $, which is dominated 
in absolute value by some integrable function with respect to the product measure $ \mu \times \mu $, and converges pointwise to the 
indicator function of the diagonal of $ X \times X $. By dominated convergence, 
\[
\lim_{t \sra 0^{+}}\int_{X \times X} F_t(x,y) \, d\mu(x) d\mu(y) = \sum_{x \in X} |\mu(\mas{x})|^2.
\]
Suppose $ ( X , g ) $ is a compact Riemannian manifold. Let $ \Delta $ be the Laplacian associated to the Riemann metric $ g $. It is a well-known 
fact \cite{Ch} that $ \Delta $ is a positive self-adjoint operator on $ L^2(X) $, with respect to the Riemannian volume measure. We let 
$ \mas{\lambda_k}_{k \geq 1} $ denote the eigenvalues of $ \Delta $, and $ \mas{\phib_k}_{k \geq 1} $ the corresponding eigenfunctions, which are orthonormal in $L^2(X)$. Note that 
the sequence of eigenvalues is non-negative and unbounded. The heat kernel $ K_t $ is the unique solution to the equation 
\[
\partial_t K_t(x,y) + \Delta_{x} K_t(x,y) = 0,
\]
with the initial condition $ \lim_{t \sra 0^{+} } K_t(x,y) = \delta_y(x) $. The heat kernel can be expressed as ( see e.g. chapter 9 in \cite{Be} )
\[
K_t(x,y) = \sum_{k \geq 1} e^{-\lambda_k t} \phib_k(x) \overline{\phib_k(y)}. 
\]
If $ ( X , g ) $ is a compact homogenous Riemannian manifold, i.e. of the form $ G / K $, where $ G $ is a Lie group, and $ K $ a closed subgroup of $ G $, the Laplacian is
equivariant under the action of $ G $ on $ L^2(X) $, and thus, $ K_t(x,x) $ is independent of $ x \in X $. If we assume that the volume has been normalized to 
be $ 1 $, we conclude that 
\[
K_t(x,x) = \sum_{k \geq 1} e^{-\lambda_k t}, 
\]
for all $ x \in X $, since the eigenfunctions are orthonormal in $ L^2(X) $. We 
now let
\[
F_t(x,y) = \frac{K_t(x,y)}{\displaystyle \sum_{k \geq 1} e^{-\lambda_k t} }.
\]
From the above discussion, it should be clear that $ F_t $ satisfies the properties above. 

If $ \mu $ is a probability measure on $ X $, we define
\[
\mu_k = \int_{X} \phib_k(x) \, d\mu(x), \quad k \geq 1.
\]
We can now formulate our first result:
\begin{lem} \label{lem:main}
Let $ X $ be a compact Riemannian homogenous space. If $ \mu $ is a probability measure on $ X $, then
\[
\lim_{t \sra \infty} \frac{\displaystyle \sum_{k \geq 1} |\mu_k|^2 e^{-\lambda_k t}}{\displaystyle \sum_{k \geq 1} e^{-\lambda_k t}} = \sum_{x \in X} |\mu(\mas{x})|^2.
\]
\end{lem}

In particular, if $ X $ is the standard $d$--dimensional torus, the lemma translates to
\[
\lim_{t \sra \infty} \frac{\displaystyle \sum_{k \in \Z^d} |\hat{\mu}(k)|^2 e^{-4\pi^2||k||^2 t}}{\displaystyle \sum_{k \in \Z^d} e^{-4\pi^2 ||k||^2 t}} = \sum_{x \in \T^d} |\mu(\mas{x})|^2,
\]
where $ \hat{\mu} $ is the Fourier transform of $ \mu $ and $ ||k|| $ is the $ l^2 $--norm of $ k \in \Z^d $. We will generalize this formula in two directions: To symmetric spaces of
compact type and to compact Heisenberg manifolds. 

\subsection{Continuous Measures on Compact Lie Groups}

Let $ G $ be a compact, simply connected and semisimple Lie group. It is a well-known fact \cite{Ta} that the set of irreducible unitary representations of $ G $ is parameterized by dominant
weights. We let $ D $ denote the set of dominant weights of the Lie algebra of $ G $, and if $ \lambda $ is in $ D $, we let $ \pi_\lambda $ denote the associated representation and
$ d_\lambda $ the dimension of this representation. Recall that all irreducible representations are finite--dimensional. The matrix coefficients of $ \pi_\lambda $ are common 
eigenvectors of the Casimir element restricted to the representation, and they all have the same eigenvalue, which we will denote by $ c_\lambda $. The character $ \chi_\lambda $
is defined to be the normalized trace of $ \pi_\lambda $, and thus $ \chi_\lambda(e) =  1 $, where $ e $ is the identity in $ G $. The heat kernel on $ G $ can be expressed as
\cite{Fe78}
\[
K_t(x,y) = \sum_{\lambda \in D} \chi_\lambda(x) \overline{\chi_\lambda(y)} \, e^{-t c_\lambda}.
\]
If $ \mu $ is a probability measure on $ G $, we define 
\[
\mu_\lambda = \int_{X} \chi_\lambda(x) \, d\mu(x).
\]
By lemma \ref{lem:main}, we now have
\begin{thm} \label{thm:cpt}
Let $ G $ be a compact, simply connected and semisimple Lie group, and suppose $ \mu $ is a probability measure on $ G $. Then,
\[
\lim_{t \sra 0^{+}} \frac{\displaystyle \sum_{\lambda \in  D} |\mu_\lambda|^2 e^{-t c_\lambda }}
{\displaystyle \sum_{\lambda \in D}  e^{-tc_\lambda}} = 
\sum_{g \in G} |\mu(\mas{g})|^2.
\]
\end{thm}

This formula should be compared to Anoussis and Bisbas result in \cite{AnBi00}. They prove, if $ G $ is a compact semi--simple Lie group and $ \mas{A_n}_{n \geq 1} $ is a sequence of subsets of dominant weights,
satisfying some technical assumptions, then
\[
\lim_{n \sra \infty} \frac{1}{|A_n|} \sum_{\lambda \in A_n} d_\lambda^{-1} ||\hat{\mu}(\pi_\lambda)||_2^2 = \sum_{x \in G} |\mu(\mas{x})|^2,
\]
where $ || \cdot || $ denotes the Hilbert--Schmidt norm, and
\[
\hat{\mu}(\pi_\lambda) = \int_{G} \pi_{\lambda}(x^{-1}) \, d\mu(x)
\]
is the operator--valued Fourier transform of $ \mu $. 

\begin{exem} 
We now consider the group $ G = SU(2) $. In this case, the dominant weights are positive half--integers \cite{Ta}, and $ c_\lambda = \frac{1}{2} \lambda ( 1 + \lambda ) $. The character associated
to the weight $ \lambda $ ( restricted to the maximal torus $ \T^1 $ ) is given by:
\[
\chi_\lambda(t) = \frac{\sin{((2\lambda+1) \pi t})}{(2\lambda+1)\sin{(\pi t)}}, \quad t \in \R,
\]
Thus, theorem \ref{thm:cpt} takes the form:
\[
\lim_{t \sra 0^{+}} \frac{\displaystyle \sum_{2\lambda \in \Z_{+}} |\mu_{\lambda}|^2 e^{-\frac{1}{2}\lambda(1+\lambda)t}}{\displaystyle \sum_{2\lambda \in \Z_{+}} e^{-\frac{1}{2}\lambda(1+\lambda)t}} = 
\sum_{x \in G} |\mu(\mas{x})|^2.
\]
\end{exem}
The extension of the above result to symmetric spaces of compact type is straight--forward.

\subsection{Continuous Measures on Nilmanifolds}
The real Heisenberg group $ H_n $ is defined to be the group 
\[
H_n = 
\left\{ 
\left(
\begin{array}{ccccc}
1 & x_1 & \ldots & x_n & z \\
 & \ddots & & & y_1 \\ 
 & & \ddots & 0 & \vdots \\
 & 0 & & \ddots & y_n \\
& & & & 1 
\end{array}
\right) 
\: | \: x,y \in \R^n, z \in \R 
\right\},
\]
and we let $ \Gamma_n $ denote the restriction of $ H_n $ to integer points. It is a well-known fact that the quotient $ X_n = H_n / \Gamma_n $ 
is a compact Riemannian manifold ( the Riemann structure is induced from Euclidean structure on the Lie algebra of $ H_n $ ). The 
spectrum of the Laplacian on $ X_n $ was calculated by Deninger and Singhof in \cite{DeSi84}. We will briefly recall their result. The spectrum is decomposed into
two parts: the first part is parameterized by $ k , h \in \Z^n $ with eigenfunctions
\[
f_{k,h}(x,y,z) = e^{2 \pi i( \l k , x \r + \l h , y \r )}
\]
with eigenvalues
\[
\lambda_{k,h} = -4\pi^2(||k||^2 + ||h||^2).
\]
The second part of the spectrum consists of eigenfunctions of the form,
\[
g_{q,m,h}(x,y,z) = e^{2 \pi i(mz + \l q , y \r)} 
\prod_{j=1}^{n} \hak{ \sum_{k \in \Z} F_{h_j}\pan{ \sqrt{2\pi |m|} \pan{x_j + \frac{q_j}{m}+k}} e^{2\pi i k m y_j}}
\]
where $ m \in \Z - \mas{0}, q \in \Z^n $ and $ h = ( h_1 , \ldots , h_n ) \in \N_0^n $, with eigenvalues
\[
\lambda_{q,m,h} = -2\pi|m|(2h_1 + \ldots + 2h_n + n + 2\pi |m|).
\]
Here, $ \mas{F_\nu}_{\nu \geq 0} $ denotes the Hermite functions,
\[
F_\nu(t) = (-1)^{\nu} e^{t^2/2} \frac{d^\nu}{dt^\nu}e^{-t^2}, \quad \nu \geq 0.
\]
As before, if $ \mu $ is a probability measure on $ X_n $, we define
\[
\mu_{k,h} = \int_{X_n} f_{k,h} d\mu
\]
and
\[
\mu_{q,m,h} = \int_{X_n} g_{k,m,h} \, d\mu.
\]
We can now formulate a generalization of Wiener's lemma on the nilmanifolds $ X_n $.
\begin{thm}
Suppose $ \mu $ is a probability measure on $ X_n $, then
\[
\lim_{t \sra 0^{+}} \frac{\displaystyle \sum_{k,h} |\mu_{k,h}|^2 e^{-t\lambda_{k,h}} + \sum_{q,m,h} |\mu_{q,m,h}|^2 e^{-t\lambda_{q,m,h}} }
{\displaystyle \sum_{k,h} e^{-t\lambda_{k,h}} + \sum_{q,m,h} e^{-t\lambda_{q,m,h}}} = \sum_{x \in X_n} |\mu(\mas{x})|^2.
\]
\end{thm}

\begin{rema}
Lemma $1$ can of course also be applied to compact locally symmetric spaces $X$, where the covering space is a symmetric space of non--compact type.
It is easy to see that the heat kernel on $X$ is the periodization of the heat kernel on the covering space. Explicit formulae for heat kernels on symmetric 
spaces of non--compact type can be found in \cite{An03}.
\end{rema}

\end{document}